\documentclass[journal]{IEEEtran} 

\usepackage{cite}

\ifCLASSINFOpdf
  \usepackage[pdftex]{graphicx}

  \DeclareGraphicsExtensions{.pdf,.jpeg,.png}
\else
   \usepackage[dvips]{graphicx}

   \DeclareGraphicsExtensions{.eps}
\fi

\usepackage[cmex10]{amsmath}

\usepackage{array}

\usepackage[caption=false,font=footnotesize]{subfig}

\usepackage{pgfplots}
\usepackage{amsfonts}
\usepackage{mathabx} 
\usepackage{amssymb}
\usepackage{enumerate}
\usepackage{latexsym}

\usepackage{url}
\usepackage{floatrow}

\usepackage[cmex10]{amsmath}
\usepackage{multirow}

\usepackage{dsfont}

\usepackage{algorithm}

\usepackage{hhline}
\newtheorem{definition}{Definition}
\newtheorem{theorem}{Theorem}

\newcommand{\tens}[1]{\boldsymbol{\mathcal{#1}}}

\newcommand{\matr}[1]{\boldsymbol{#1}}

\newcommand{\vect}[1]{\boldsymbol{#1}}

\newcommand{\tensp}{\mathop{\otimes}}      
\newcommand{\kr}{\odot}     
\newcommand{\kron}{\mathop{\boxtimes}}   
\newcommand{\pout}{\circ}   

\newcommand{\vecl}{\text{vec}}

\usepackage{color}

\definecolor{burgundy}{RGB}{159,29,53}
\definecolor{arylideyellow}{rgb}{0.91, 0.84, 0.42}
\definecolor{bananayellow}{rgb}{1.0, 0.88, 0.21}

\makeatletter
\renewcommand{\maketag@@@}[1]{\hbox{\m@th\normalsize\normalfont#1}}%
\makeatother

\begin{document}

\title{About Notations in Multiway Array Processing}

\author{J\'er\'emy~E.~Cohen} 


\maketitle

\begin{abstract}
This paper gives an overview of notations used in multiway array processing. We redefine the vectorization and matricization operators to comply with some properties of the Kronecker product. The tensor product and Kronecker product are also represented with two different symbols, and it is shown how these notations lead to clearer expressions for multiway array operations. Finally, the paper recalls the useful yet widely unknown properties of the array normal law with suggested notations.
\end{abstract}


%

\section*{Introduction}
Since tensors have become a popular topic in data science and signal processing, a large number of papers have been published to describe at best their different properties \cite{Kolda2009,kolda2006multilinear, DeLathawer2000, cichocki2009nonnegative, Bro1998, kiers2000towards, kroonenberg1983three}. Yet a consensus on notations has not been found. Some authors refer to the Tucker and Kruskal operators \cite{kruskal1977three, kolda2006multilinear} while others never use these notations but make use of the n-way product \cite{DeLathawer2000, hoff2011separable}. There is at least three different unfolding methods in the literature \cite{Kolda2009, DeLathawer2000, landsberg2012tensors}. The Kronecker product alone has two notations in the community \cite{Comon2014,landsberg2012tensors, Kolda2009}.

The tensor product itself is almost never used, even though it is the very foundation of tensor algebra. Instead authors sometimes refer to the outer product $\pout$, which may be seen as a tensor product of vectors in the canonical basis of each vector space. But the tensor product allows for wider generalization since it applies to tensors of any order and without referring to any basis. Some authors suggest the use of two different symbols for the Kronecker product $\kron$ and the general tensor product $\tensp$ \cite{Comon2014}. But the same symbol has been used historically and through the literature \cite{brewer1978kronecker}.

The Kronecker product is the expression of the tensor product for matrices when a basis has been given. Yet using the same symbol leads to confusion when manipulating arrays and multilinear operators at the same time. Moreover, manipulating the arrays with matricizations and vectorizations exactly means that the difference between a general tensor product and a basis-dependent Kronecker product is of crucial importance. This leads also to redefining some well known operations, namely the matricization and the vectorization.

The main goal of this paper is therefore to set notations for array and manipulations on arrays that are consistent with one another, and with notations used in quantum physics and algebraic geometry where tensors have been used for decades  \cite{nielsen2010quantum}. We recommend that the vectorization and matricization operators are computed so that the tensor products and Kronecker products do not have to be swapped versions of each others, which is compatible with definitions from algebraic geometry \cite{landsberg2012tensors}.  In other words, the suggested notations will lead to \[ \vecl\left(\vect{a}\tensp{\vect{b}}\right)=\vect{a}\kron\vect{b},\] where $\tensp$ is the outer product, and $\kron$ is the Kronecker product. This differs from the usual equality $ \vecl\left(\vect{a}\tensp{\vect{b}}\right)=\vect{b}\kron\vect{a}$, which leads to unnecessarily complicated equations.

The first section introduces some notations for all operators used extensively in multiway array processing. Then some useful properties of tensors are given in Section \ref{section2} using the suggested notations. Finally Section \ref{section3} exposes the array normal law as defined by Hoff \cite{hoff2011separable}. It is shown that the array normal law can be a handy tool for multiway array analysis.


\section{A Few Definitions}
\label{section1}
 \subsection{About tensor products}
First let us define the tensor product as given in \cite{bourbaki1998algebra}. 

\begin{definition}
Let $\mathcal{E}$, $\mathcal{F}$ be two vector spaces on a field $\mathcal{K}$. There is a vector space $\mathcal{E} \tensp \mathcal{F}$, called the tensor product space, and a bilinear mapping 
\[
\tensp ~ : ~ \mathcal{E} \times \mathcal{F} \rightarrow \mathcal{E} \tensp \mathcal{F} 
\] so that for every vector space $\mathcal{G}$ and for all bilinear mapping $g$ from $\mathcal{E}\times \mathcal{F}$ to $\mathcal{G}$, there is one and only one linear mapping $h$ from $\mathcal{E}\tensp \mathcal{F}$ to $\mathcal{G}$ defined by $g(x,y)=h(x\tensp y)$. This extends to multilinear mappings.
\end{definition} In other words, the tensor product of vector spaces $\left(\mathcal{E}_i\right)_{i\leq N}$ builds one (linear) vector space $\bigotimes_i \mathcal{E}_i$.

This definition is not basis dependent, which means that the tensor product applies not solely to arrays, and that real-valued tensors are items from a tensor space $\bigotimes_{i\leq N} \mathds{R}^{n_i}$. When a basis is given for the tensor space, this represents a tensor by an array of coordinates. Two well known basis-dependent representation for the tensor product are the outer product and the Kronecker product \cite{DeLathawer2000,brewer1978kronecker}. First let us define the Kronecker product of two arrays :
\begin{definition}
The Kronecker product of two arrays $\matr{A} \in \mathds{R}^{p_1\times q_1}$ and $\matr{B} \in \mathds{R}^{p_2\times q_2}$ is denoted by $\matr{A}\kron\matr{B} \in \mathds{R}^{p_1p_2 \times q_1q_2}$ and is defined by:
\begin{equation*}
\matr{A} \kron \matr{B} := \left[
\begin{array}{c|c|c|c}
a_{11}\matr{B} & a_{12}\matr{B} & \dots & a_{1q_1}\matr{B} \\
\hline
a_{21}\matr{B} &  &  & \\
\hline
\vdots & & & \\
\hline
a_{p_11}\matr{B} & &  & a_{p_1q_1}\matr{B} 
\end{array} \right].
\end{equation*}
\end{definition}
We also give the definition of the outer product of two vectors in a given basis.
\begin{definition}
The outer product of two vectors $\vect{a} \in \mathds{R}^{p_1}$ and $\vect{b} \in \mathds{R}^{p_2}$ is denoted by $\vect{a}\tensp\vect{b} \in \mathds{R}^{p_1\times p_2}$ and is defined by :
\begin{equation*}
\left(\vect{a} \tensp \vect{b}\right)_{ij} = a_ib_j
\end{equation*}
\end{definition}
These definitions can be extended in a trivial manner to more than two arrays.

Ambiguity happens when the used tensor product maps to the matrix algebra, but no basis is given. We then need to use the symbol $\tensp$ for this tensor product. Moreover, the Kronecker product and the tensor product in matrix algebra coincide when a basis is provided. It is shown in the rest of this paper that this notation eases array manipulation significantly. 

On the other hand, the outer product does not need a third symbol, since there is no ambiguity with the tensor product. Whenever dealing with arrays, the tensor product of vectors is to be read as an outer product, but in any other case the outer product makes no sense. To comply with usual notations in multilinear algebra, the symbol $\pout$ both proposed and later abandoned by Comon \cite{comon1998blind,stegeman2010subtracting} should be dropped since it refers also to composition. Also, the usual notation for the tensor product is $\tensp$, except in part of the multiway array processing community.

Multiway array processing is based on applying multilinear operators, which are themselves higher order tensors (see appendix \ref{app1}), on data arrays to mine these data. A well-known example is the multilinear operator obtained by High Order Singular Value Decomposition (HOSVD) \cite{DeLathawer2000} yielding compression for the data array. However, when fixing a basis there is not a single definition of a product of two tensors. In particular, if $\matr{U}_i$ are matrices in $\mathds{R}^{q_i\times n_i}$, the operator  $\bigotimes_{i\leq N}\matr{U}_i$ acting on tensors in $\bigotimes_{i\leq N} \mathds{R}^{n_i}$ cannot be expressed in a given basis with outer products and Kronecker products \cite{lekheng2015}. Still for computing this product of tensors, some authors defined the $k$-way product  \cite{kolda2006multilinear,DeLathawer2000}  :
\begin{equation*}
\left(\tens{Y}\bullet_k\matr{U}_k\right)_{i_1, \dots, i_N} := \sum_{j=1}^{n_k}{Y_{i_1, \dots, i_{k-1}, j, i_{k+1}, \dots, i_N}\matr{U}_{i_n, j}}
\end{equation*}
The bullet symbol $\bullet$ should be preferred to the multiplication symbol $\times$, since the latter may already have other meanings. Moreover, the $k$-way product is the contraction of a matrix along its second mode with a tensor along its $k$-th mode.

Now operators $\bigotimes_{i\leq N}\matr{U}_i$ are already extensively used in quantum mechanics among others, and their action on a tensor is denoted with an implicit multilinear product:
\begin{equation*}
\left( \bigotimes_{i=1}^N\matr{U}_i \right) \tens{Y} := \tens{Y} \bullet_1  \matr{U}_1  \dots \bullet_N  \matr{U}_N
\end{equation*}
which is to be read as $N$ contractions for the second mode of each matrix in the tensor product with the $i$-th mode of the tensor. This also provides a computation mean for the application of an operator on multiway arrays. One can simply sequentially unfold the tensor in each mode, multiply the resulting matrix with the matrix operator $\matr{U}_i$ expressed in the right basis, and refold the tensor. For matrices, the 2-mode product is indeed well known :
\begin{equation*}
\matr{M}\bullet_1\matr{U}_1\bullet_2\matr{U}_2=\matr{U}_1\matr{M}\matr{U}_2^T
\end{equation*}
One last operation useful for multiway array processing is the (basis-dependent) Khatri-Rao product :
\begin{equation*}
\matr{A} \kr \matr{B} = \left[\vect{A}_{:1} \kron \vect{B}_{:1}, \dots, \vect{A}_{:p} \kron \vect{B}_{:p} \right]
\end{equation*}

	\subsection{About manipulation of arrays}

A major claim of this paper is to modify the usual definition of the vectorization operator for arrays \cite{neudecker1969some}. For now, the usual definition states that an array should be vectorized choosing elements by the reverse lexicographic order of indices, \textit{i.e.} columns after columns in the first mode, sliding along the second mode, then the third mode and so on (see figure \ref{fig1}a). 

Yet every tensor can be expressed as a sum of $R$ tensor products of vectors by the following model, called Canonical Polyadic (CP) decomposition or PARAFAC \cite{carroll1970analysis,harshman1970foundations} :
\begin{equation*}
\tens{Y} = \sum\limits_{r=1}^{R}{\bigotimes_{i=1}^{N} \vect{a}_r^{(i)}}
\end{equation*}
and using the columnwise vectorization \textit{is not compatible} with the way the Kronecker product and the CP decomposition are defined. To see this, take a matrix $\matr{M}$ defined by an outer product $\matr{M}=\vect{a} \tensp \vect{b}$. When vectorizing $\matr{M}$ columnwise, it is well known that vectors $\vect{a}$ and $\vect{b}$ are swapped to express $\matr{M}$ as a Kronecker product, since in its definition the Kronecker product makes the second index vary first : $\text{vec}(\matr{M})=\vect{b}\kron\vect{a}$. The point here is that once the Kronecker product has been defined, it fixes the direction for the vectorization of any array in order to get a direct property like $\text{vec}\left(\bigotimes_i \vect{a}_i\right)=\bigboxtimes_i \vect{a}_i$.

Thus we give the following definition of the vectorization operator of an array, along with some basic properties in Table \ref{table1}:
\begin{definition}
The vectorization of an array $\tens{Y}$, denoted $\vecl(\tens{Y})$, is the isomorphism defined as the following function of its CP decomposition:
\begin{align*}
   \mathds{R}^{n_1 \times \dots \times n_N}   & \rightarrow \mathds{R}^{\prod\limits_{i=1}^{N}{n_i} \times 1} \\
 \text{vec} : \quad  \sum\limits_{r=1}^{R}{\bigotimes_{i=1}^{N} \vect{a}_r^{(i)}} & \mapsto   \sum\limits_{r=1}^{R}{\bigboxtimes_{i=1}^{N} \vect{a}_r^{(i)}}
\end{align*} 
This definition means the elements in the array are taken in the lexicographic order of their index, \textit{i.e.} along the last mode, then along the mode $N-1$ and so on (see figure \ref{fig1}b).
\end{definition}

In the same spirit, we define the matricization (or the unfolding) of an array similarly to what is called flattening in algebraic geometry \cite{landsberg2012tensors}. Some properties are given in table \ref{table1}.
\begin{definition}
The matricization of an array $\tens{Y}$ along the $i$-th mode, denoted $\left[ \tens{Y}\right]_{(i)}$ is the isomorphism defined as the following function of its CP decomposition :
\begin{align*}
   \mathds{R}^{n_1 \times \dots \times n_N}   & \rightarrow \mathds{R}^{n_i \times \prod\limits_{j\neq i}^{N}{n_j}} \\
 \left[ . \right]_{(i)} : \quad  \sum\limits_{r=1}^{R}{\bigotimes_{j=1}^{N} \vect{a}_r^{(j)}} & \mapsto   \sum\limits_{r=1}^{R}{\vect{a}_r^{(i)} \tensp \bigboxtimes_{j\neq i}^{N} \vect{a}_r^{(j)}}
\end{align*} 
\end{definition}

\begin{figure}
\begin{tikzpicture}[scale=1.6]

\draw[very thick,burgundy] (0,0,0)--(0,-1.05,0);
\draw[->,thick,black] (-0.1,-0.3,0)--(-0.1,-0.7,0);
\node at (-0.3,-0.5,0) {(1)};
\node at (0,-1.25,0) {i};

\draw[very thick,burgundy] (0,0,0)--(1,0,0);
\draw[->,thick,black] (0.5,0.1,0)--(0.85,0.1,0);
\node at (0.7,0.3,0) {(2)};
\node at (1.15,0,0) {j};

\draw[very thick,burgundy] (0,0,0)--(0,0,-0.7);
\draw[->,thick,black] (-0.15,0,-0.3)--(-0.15,0,-0.7);
\node at (-0.4,0,-0.6) {(3)};
\node at (-0.1,0,-1.1) {k};

\node at  (0.5,-1.9,0) {\footnotesize a) columnwise vectorization};


\draw[very thick,burgundy] (3,0,0)--(3,-1.05,0);
\draw[->,thick,black] (2.9,-0.3,0)--(2.9,-0.7,0);
\node at (2.7,-0.5,0) {(3)};
\node at (3,-1.25,0) {i};

\draw[very thick,burgundy] (3,0,0)--(4,0,0);
\draw[->,thick,black] (3.5,0.1,0)--(3.85,0.1,0);
\node at (3.7,0.3,0) {(2)};
\node at (4.15,0,0) {j};

\draw[very thick,burgundy] (3,0,0)--(3,0,-0.7);
\draw[->,thick,black] (2.85,0,-0.3)--(2.85,0,-0.7);
\node at (2.6,0,-0.6) {(1)};
\node at (2.9,0,-1.1) {k};

\node at  (3.5,-1.9,0) {\footnotesize b) suggested vectorization};

\end{tikzpicture}

\caption{Two possibilities for vectorizing an array. If $Y_{ijk}$ is an element in array $\tens{Y}$ of size $I\times J\times K$, then in the vectorized form $\vect{y}$ it is located at : a) $y_{(k-1)IJ+(j-1)I+i}$ \quad b) $y_{(i-1)JK+(j-1)K+k}$}
\label{fig1}
\end{figure}
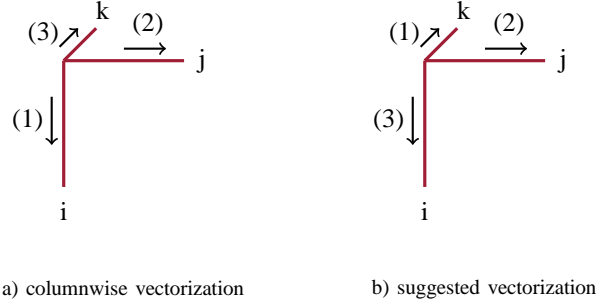

\begin{figure}
\begin{tikzpicture}[scale=1.6]

\draw[fill=burgundy] (0,0,0)--(0,-1.05,0)--(0.1,-1.05,0)-- (0.1,-1.05,-0.7)--(0.1,0,-0.7)--(0,0,-0.7)--(0,0,0);

\draw (0.1,0,0)--(0.1,-1.05,0);
\draw (0.1,0,0)--(0.1,0,-0.7);

\node at (0.1,-1.2,0) {\tiny$\matr{M}_1$};

\draw[fill=burgundy] (0.9,0,0)--(0.9,-1.05,0)--(1,-1.05,0)-- (1,-1.05,-0.7)--(1,0,-0.7)--(0.9,0,-0.7)--(0.9,0,0);
\draw (1,0,0)--(1,-1.05,0);
\draw (1,0,0)--(1,0,-0.7);
\node at (1,-1.2,0) {\tiny$\matr{M}_K$};

\node at (0.6,-0.45,0) {\dots};

\draw[very thick] (0,0,0)--(0,-1.05,0);

\draw[very thick] (0,0,0)--(1,0,0);

\draw[very thick] (0,0,0)--(0,0,-0.7);
\node at  (2,-1.9,0) {$
\begin{aligned}
\matr{Y}_{(1)}&=&\left[ \color{burgundy}\matr{M}_1 \color{black}| \dots | \color{burgundy}\matr{M}_K\color{black} \right] \\
\matr{Y}_{(2)}&=&\left[ \color{arylideyellow!90!black}\matr{N}_1 \color{black}| \dots | \color{arylideyellow!90!black}\matr{N}_I\color{black} \right]  \\
\matr{Y}_{(3)}&=&[ \color{arylideyellow!90!black}\matr{N}_1^T \color{black}| \dots | \color{arylideyellow!90!black}\matr{N}_I^T\color{black} ]  \end{aligned}$};
\node at  (2,-3,0) {where $\matr{M}_{i} \in \mathds{R}^{I\times K} $ and $\matr{N}_{i} \in \mathds{R}^{J\times K}$};

\draw[fill=arylideyellow] (3,0,0)--(3,-0.1,0)--(4,-0.1,0)-- (4,-0.1,-0.7)--(4,0,-0.7)--(3,0,-0.7)--(3,0,0);

\draw (4,0,0)--(4,-0.1,0);
\draw (4,0,0)--(4,0,-0.7);

\node at (4.45,0,-0.2) {\tiny$\matr{N}_1$};

\draw[fill=arylideyellow] (3,-0.95,0)--(3,-1.05,0)--(4,-1.05,0)-- (4,-1.05,-0.7)--(4,-0.95,-0.7)--(3,-0.95,-0.7)--(3,-0.95,0);
\draw (4,-0.95,0)--(4,-1.05,0);
\draw (4,-0.95,0)--(4,-0.95,-0.7);
\draw (3,-0.95,0)--(4,-0.95,0);
\node at (4.45,-0.95,-0.2) {\tiny$\matr{N}_I$};

\node at (3.6,-0.35,0) {\vdots};

\draw[very thick] (3,0,0)--(3,-1.05,0);

\draw[very thick] (3,0,0)--(4,0,0);

\draw[very thick] (3,0,0)--(3,0,-0.7);


\end{tikzpicture}
\caption{Unfoldings of a three-way array $\tens{Y}$ of size $I\times J\times K$} 
\label{figure2}
\end{figure}

This definition can be exploited to obtain the unfoldings of a three way array $\tens{Y}$ as shown in figure \ref{figure2}. For example, a compact MATLAB R2014b code for computing the three unfoldings is:
\begin{align*}
\matr{Y}_{(1)} &=& \text{reshape}\left(\text{permute} \left(\tens{Y},[1,3,2] \right),n_1,n_2n_3 \right) \\ 
\matr{Y}_{(2)} &=& \text{reshape}\left(\text{permute} \left(\tens{Y},[2,3,1] \right),n_2,n_1n_3 \right) \\ 
\matr{Y}_{(3)} &=& \text{reshape}\left(\text{permute} \left(\tens{Y},[3,2,1] \right),n_3,n_1n_2 \right) 
\end{align*}

\section{Some Useful Properties}
\label{section2}

When manipulating multiway arrays, there are multiple properties that may simplify calculations by wide margins. Formula (\ref{diff}) is useful for obtaining a compact formula for the Jacobian of a CP model or obtaining bounds on estimation errors. Table \ref{table1} provides overall well known results \cite{brewer1978kronecker}, with our suggested notations. Note that permutation of vectors or operators when switching from tensor products to vectorized arrays do not appear anymore.

In what follows, some well-studied decompositions are revisited with the notations we promote.

\paragraph{Exemple 1} The Tucker decomposition of a tensor $\tens{Y}$
\[
Y_{i_1 \dots i_N }=\sum\limits_{r_1 \dots r_N}^{R_1 \dots R_N}{\left( \prod\limits_{k=1}^{N}{U_{i_k r_k}}\right) G_{r_1 \dots r_N} }
\]
can be expressed independently of any fixed basis as the following :

\[
\tens{Y} = \left( \bigotimes_{i=1}^N \matr{U}_i \right) \tens{G}
\]
where operators $\matr{U}_i$ may be supposed unitary for identifiability of the model.

Thus finding the Tucker decomposition of a tensor means finding linear operators $\matr{U}_i$ in $\mathds{R}^{n_i\times R_i}$ so that the multilinear operator $\left( \bigotimes_{i=1}^N \matr{U}^T_i \right)$ projects the tensor on its true subspace $\bigotimes_i \mathds{R}^{R_i}$, in the same spirit as Principal Component Analysis.

Also from Table \ref{table1}, unfoldings can be easily obtained:

\begin{align*}
\matr{Y}_{(i)} &= \left( \matr{U}_i \otimes \bigotimes_{j\neq i}^N \matr{U}_j \right) \matr{G}_{(i)} \\
                     &= \matr{U}_i \matr{G}_{(i)} \left( \bigboxtimes_{j\neq i}^N \matr{U}_j^T \right) \\
\end{align*}

Of course for actual computation of the action of $\left( \bigotimes_{i=1}^N \matr{U}_i^T \right) $ on $\tens{Y}$, all the $i$-way $\tens{Y} \bullet_i \matr{U}_i^T$ products may first be computed sequentially before matricizing the tensor. 

\paragraph{Example 2} In the first section, the CP decomposition of a tensor was defined. With our notations, this model can be easily expressed as follows:
\[
\tens{Y} = \sum\limits_{r=1}^{R}{\bigotimes_{i=1}^{N} \vect{a}_r^{(i)}} = \left( \bigotimes_{i=1}^N \matr{A}_i \right) \tens{I}
\]
where $\matr{A}_i = \left[ \vect{a}_1 \dots \vect{a}_R \right]$ and $\tens{I}$ is the diagonal tensor with entries equal to $1$.  The CP model can be described with unfoldings and in a vectorized form using (\ref{CPmodel}) from table \ref{table1} :

\begin{align*}
\matr{Y}_{(i)} & = \matr{A}_i  \tens{I}_{(i)} \left( \bigboxtimes_{j \neq i}^N \matr{A}_j \right)^T =\matr{A}_i \left( \bigodot_{j \neq i}^N \matr{A}_j \right)^T \\
\vecl(\tens{Y}) & = \sum\limits_{r=1}^{R}{\bigboxtimes_{i=1}^N \vect{a}_r^{(i)}} = \left( \bigodot_{i=1}^N \matr{A}_i \right) \mathds{1}
\end{align*}
where $\mathds{1}$ is a vector of ones.

\begin{table}
\begin{align}
& \vecl (\vect{a} \tensp \vect{b})  = \vect{a} \kron \vect{b} \\
& \vecl \left( \matr{A}\matr{X}\matr{B}^T\right)  = \left(\vect{A} \kron \vect{B}\right) \vecl \left(\matr{X}\right) \\
& \left( \matr{A} \tensp \matr{B} \right) \left( \matr{C} \tensp \matr{D} \right) = \left( \matr{A}\matr{C} \tensp \matr{B}\matr{D} \right), \text{ if compatible} \\
& \left( \bigotimes_{i=1}^{N} \matr{U}_i  \right) \left( \bigotimes_{i=1}^{N} \vect{a}^{(i)} \right) = \left( \bigotimes_{i=1}^{N} \matr{U}_i\vect{a}^{(i)} \right) \\
& \left( \matr{A} \kron \matr{B} \right) \left( \matr{C} \odot \matr{D} \right) = \left( \matr{A}\matr{C} \odot \matr{B}\matr{D} \right), \text{ if compatible} \\
& \left( \matr{A} \tensp \matr{B} \right)^T =  \matr{A}^T \tensp \matr{B}^T \\
& \matr{A}\matr{X} + \matr{XB} = \matr{D} \Leftrightarrow \left( \matr{A} \kron \matr{B} \right) \vecl\left(\matr{X}\right) = \vecl\left(\matr{D}\right) \\
& \vecl\left( \bigotimes_{i=1}^N \vect{a}^{(i)} \right) = \bigboxtimes_{i=1}^N \vect{a}^{(i)} \text{   in the same order}\\
& \tens{Y} = \tens{X} \bullet_i \matr{U}  \Leftrightarrow \matr{Y}_{(i)}=\matr{U}\matr{X}_{(i)} \\ 
& \vecl \left(\left( \bigotimes_{i=1}^N \matr{U}_i \right)  \tens{Y} \right)=  \left( \bigboxtimes_{i=1}^N \matr{U}_i \right) \vecl\left( \tens{Y} \right) \\
& \left[ \left( \bigotimes_{j=1}^N \matr{U}_j \right) \tens{Y} \right]_{(i)} = \matr{U}_i\matr{Y}_{(i)}\left(\bigboxtimes_{j\neq i}^N \matr{U}_j^T \right) \\
& \left[ \left( \bigotimes_{j=1}^N \matr{U}_j \right) \tens{Y} \right]_{(i)} = \left(\matr{U}_i \tensp \left( \bigotimes_{j\neq i}^N \matr{U}_j^T \right) \right) \matr{Y}_{(i)} \\
& \tens{Y} = \sum\limits_{r=1}^{R}{\bigotimes_{j=1}^N \vect{a}_r^{(j)}} \Leftrightarrow \matr{Y}_{(i)} = \sum\limits_{r=1}^{R}{\vect{a}^{(i)}_r \tensp \bigboxtimes_{j\neq i}^N \vect{a}_r^{(i)}} \\
& \sum\limits_{r=1}^{R}{\vect{a}_r \kron \vect{b}_r} = \left(\matr{A} \odot \matr{B}\right)\mathds{1}, \text{ where } \matr{A}=\left[\vect{a}_1 \dots \vect{a}_R\right] \\
& \matr{Y}_{(i)} = \sum\limits_{r=1}^{R}{\vect{a}^{(i)}_r \tensp \bigboxtimes_{j\neq i}^N \vect{a}_r^{(j)}} = \matr{A}_i\left(\bigodot_{j \neq i}^N \matr{A}_j^T \right) \label{CPmodel}\\
& \frac{\partial \bigotimes_{j=1}^N \vect{a}_j}{\partial \vect{a}_i} = \vect{a}_1 \tensp \dots \tensp \vect{a}_{i-1} \tensp \matr{I}  \tensp \dots \tensp \vect{a}_N
\label{diff} \\
& \left| \bigotimes_{i=1}^N \matr{U}_i \right| = \prod_{i=1}^N \left| \matr{U}_i\right|^{\prod_{j \neq i}^N n_j} 
\end{align}
\caption{Some useful properties}
\label{table1}
\end{table}


Applying a multilinear transformation to the tensor results in operating $\bigotimes_{i=1}^N \matr{W}_i$ on $\tens{Y}$, which has the following CP model :

\[
\left( \bigotimes_{i=1}^N \matr{W}_i \right) \tens{Y} =  \sum\limits_{r=1}^{R}{\bigotimes_{i=1}^{N} \matr{W}_i\vect{a}_r^{(i)}}.
\]

This kind of preprocessing is very useful if the observation noise is correlated, or in some applications e.g. flexible coupled tensor decomposition \cite{Farias2015} where the decomposition may be computed in a transformed domain for one of the mode.


\section{Second Order Statistics for Tensors}
\label{section3}

When dealing with Gaussian non-i.i.d. multiway arrays, the most natural way of describing the distribution would be to give a definition to array normal laws. This has been done for matrices \cite{gupta1999matrix}, and recently for arrays of any order \cite{hoff2011separable}. The definition from \cite{hoff2011separable} is the following~:

\begin{definition}
Let $\tens{X}$ be a multivariate random variable in $\mathds{R}^{n_1\times \dots \times n_N}$. We say that $\tens{X}$ follows an array normal law of mean $\tens{M}$ and with tensor covariance $\matr{\Gamma}=\bigotimes_{i=1}^N \matr{\Sigma}_i$ if and only if 
\[
p\left(\tens{X}|\tens{M},\matr{\Gamma} \right) = \frac{\text{exp}\left(- \frac{\| \matr{\Gamma}^{-\frac{1}{2}}(\tens{X}-\tens{M}) \|^2}{2} \right)}{(2\pi)^{\prod_i{\frac{n_i}{2}}} \left|\matr{\Gamma}\right|^{\frac{1}{2}} }
\] where $\matr{\Gamma}^{-\frac{1}{2}} = \bigotimes_{i=1}^N \matr{\Sigma}_i^{-\frac{1}{2}}$ and $\matr{\Sigma}_i$ are symmetric. It is noted as $\tens{X} \sim \mathcal{A}\mathcal{N}\left( \tens{M}, \matr{\Gamma} \right)$
\end{definition}

Array normal distributions are useful when manipulating the data arrays before decomposing them. They suppose the covariance is expressed in every mode, \textit{i.e.} the covariance of the vectorized tensor can be expressed as a Kronecker product of $N$ symmetric matrices, \textit{i.e.} the covariance is separable. 

For example, say it is needed to preprocess some tensor $\tens{Y}$, noisy version of $\tens{X}$ corrupted by i.i.d. Gaussian noise, with the multilinear operator $\bigotimes_{i\leq N}\matr{U}_i$. Then the array $\tens{T}=\left(\bigotimes_{i\leq N}\matr{U}_i\right)\tens{Y}$ follows the array normal law 
\[ \tens{T} \sim \mathcal{AN}\left( \left(\bigotimes_{i=1}^N\matr{U}_i\right)\tens{X}, \bigotimes_{i=1}^N \matr{U}_i\matr{U}_i^T \right). \]

Now with notations from two way array processing \cite{gupta1999matrix}, the matrix law of the $i_{th}$ unfolding is given by :
\[
\matr{T}_{(i)} \sim \mathcal{MN}\left( \matr{U}_i \matr{X}_{(i)} \bigboxtimes_{j \neq i}^N \matr{U}_j^T , \matr{U}_i\matr{U}_i^T, \bigboxtimes_{j\neq i}^N \matr{U}_j\matr{U}_j^T \right)
\]
and the normal law of the vectorized tensor is :
\[
\vecl(\tens{T}) \sim \mathcal{N}\left( \left(\bigboxtimes_{i=1}^N\matr{U}_i \right) \vecl(\tens{X}) , \bigboxtimes_{i=1}^N\matr{U}_i\matr{U}_i^T \right)
\]

A direct consequence is that preprocessing one mode of a noisy tensor modifies the covariance of the noise in this mode but this does not affect the other modes. In a least square procedure for fitting the CP decomposition, this results in modifying the maximum likelihood estimate for the modified mode only.

\section*{Concluding remarks}
This paper proposed a reworked set of notations for the multiway array processing community. The notations commonly used are sometimes not compatible with other communities, and the author believes the suggested notations simplify calculations with arrays. The permutation induced by the column wise vectorization operator is a source of error, as well as the confusion between tensor products and the Kronecker product. Moreover, unfolding operators given in the literature suggest a similar permutation, which is also unnecessary. With proposed notations, there are no more permutations anywhere in the computation of matricized and vectorized tensors. Finally, this framework enables a simple definition of the multivariate normal distribution for arrays of any order, which is useful for example in computing decompositions of noisy preprocessed tensors.

Notations in this paper are of course one choice among others, and some definitions are not the ones usually used. For digging further into tensors, some nice surveys are available \cite{Kolda2009,Comon2014,grasedyck2013literature}. 






\appendices

\section{Linear operators acting on tensors}
\label{app1}

We state and prove that linear operators acting on a tensor space of finite dimension linear spaces is itself a tensor space. This justifies the notation $\bigotimes_{i=1}^N\matr{U}_i$ suggested in this paper. A deeper proof including infinite dimensions can be found in the excellent book by Hackbusch \cite{hackbusch2012tensor}.

Let $\mathcal{E}$, $\mathcal{E}'$, $\mathcal{F}$ and $\mathcal{F}'$ be four finite vector spaces on a field $\mathcal{K}$, of respective dimensions $n,n$ and $m,m$. Let $\otimes$ and $\otimes'$ be two tensor products on $\mathcal{E} \times \mathcal{F}$ and $\mathcal{E}' \times \mathcal{F}'$, and consider the following mapping:
\begin{equation*}
\begin{array}{cccc}
&   \mathcal{L}(\mathcal{E},\mathcal{E}') \times \mathcal{L}(\mathcal{F},\mathcal{F}')  & \rightarrow  &\mathcal{L}(\mathcal{E} \otimes \mathcal{F}, \mathcal{E}' \otimes' \mathcal{F}') \\
 \otimes_p : & \quad  (u,v) & \mapsto  & u\otimes_p v : (x\otimes y) \mapsto u(x) \otimes' v(y)
 \end{array}
\end{equation*} 
where $\mathcal{L}(\mathcal{E},\mathcal{E}')$ is the linear space of linear operators mapping $\mathcal{E}$ to $\mathcal{E}'$.

\begin{theorem}
$\mathcal{L}(\mathcal{E} \otimes \mathcal{F}, \mathcal{E}' \otimes' \mathcal{F}')$ associated with the bilinear mapping $\otimes_p$ is a tensor space, i.e. 
\[
\mathcal{L}(\mathcal{E} \otimes \mathcal{F}, \mathcal{E}' \otimes' \mathcal{F}') = \mathcal{L}(\mathcal{E}, \mathcal{E}') \otimes_p \mathcal{L}(\mathcal{F},\mathcal{F}').
\]
\end{theorem}

To prove theorem 1, we only need to check that $\otimes_p$ maps one basis of $\mathcal{L}(\mathcal{E},\mathcal{E}') \times \mathcal{L}(\mathcal{F},\mathcal{F}')$ to a free family of $\mathcal{L}(\mathcal{E} \otimes \mathcal{F}, \mathcal{E}' \otimes' \mathcal{F}')$, which yields injectivity. Then we will be able to conclude arguing that this two spaces have the same dimension, so that by Green's theorem, $\otimes_p$ is a bijective bilinear map from the carthesian product space to a linear space, which is exactly what a tensor product is \cite{bourbaki1998algebra}.

Let $\{e_i\}_i$, $\{e'_i\}_i$, $\{f_j\}_j$ and $\{f'_j\}_j$ be some bases of $\mathcal{E}, \mathcal{E}', \mathcal{F}$ and $\mathcal{F}'$. Define the following basis for $\mathcal{L}(\mathcal{E}, \mathcal{E}' )$:
\[
u_i(e_j) =\delta_{ij} e'_i
\]
where $\delta_{ij}$ is the kronecker symbol equal to 1 if and only if $i$ equals $j$. Define the basis $\{ v_i \}$ similarly. Let us prove that $\{ u_i \otimes_p v_j \}_{ij}$ is a free family.

Given some $\{ \lambda_{ij} \}_{ij}$, suppose for any $x \otimes y$ in $\mathcal{E} \otimes \mathcal{F}$
\[
\sum\limits_{i,j}{\lambda_{ij} u_i(x) \otimes v_j(y)} = 0 .
\]
This is equivalent to 
\[
\sum\limits_{i,j}{\lambda_{ij} \nu_i \nu'_j e'_i \otimes f'_j} = 0 .
\]
where $\nu_i$ and $\nu'_j$ are coefficients of $x$ and $y$ in bases $\{e_i\}_i$ and $\{f_j\}_j$. Since this is true for any $x,y$ and that $\{e'_i \otimes' f'_j\}_{ij}$ is a basis of $\mathcal{E}' \otimes' \mathcal{F}'$, all lambdas have to go to zero, thus proving theorem 1.
\section*{Acknowledgment}
The author would like to thank Pierre Comon and Rodrigo Cabral Farias for useful discussions.

\ifCLASSOPTIONcaptionsoff
  \newpage
\fi

\bibliographystyle{IEEEtran}
\bibliography{biblio}

\vfill


\end{document}